\documentclass[a4paper,11pt]{article}
\usepackage{amsfonts,amssymb,amscd,latexsym}
\usepackage[cp1251]{inputenc}
\usepackage[russian]{babel}
\begin{document}
\Large
\newcounter{num}[section]
\setcounter{num}{0}
\renewcommand{\thenum}{\arabic{section}.\arabic{num}}

\newcommand{\gx}{{\frak g}}
\newcommand{\hx}{{\frak h}}
\newcommand{\px}{{\frak p}}
\newcommand{\qx}{{\frak q}}
\newcommand{\Rx}{{\frak R}}
\newcommand{\rx}{{\frak r}}
\newcommand{\Ac}{{\cal A}}

\newcommand{\pv}{{\frak p}{\frak v}}
\newcommand{\qv}{{\frak q}{\frak v}}
\newcommand{\al}{\alpha}
\newcommand{\la}{\lambda}

\renewcommand{\leq}{\leqslant}
\renewcommand{\geq}{\geqslant}

\newcommand{\Irr}{{\mathrm{Irr}}}
\newcommand{\spn}{{\mathrm{span}}}
\newcommand{\ad}{{\mathrm{ad}}}
\newcommand{\Spec}{{\mathrm{Spec}}}
\newcommand{\ind}{{\mathrm{ind}}}
\newcommand{\Img}{{\mathrm{Im}}}
\newcommand{\Ker}{{\mathrm{Ker}}}
 \newcommand{\tr}{{\mathrm{tr}}}
 \newcommand{\Prim}{{\mathrm{Prim}}}
 \newcommand{\Dix}{{\mathrm{Dix}}}
 \newcommand{\KKer}{{\frak K}{\frak e}{\frak r}}

\newcommand{\Zb}{{\Bbb Z}}

\newcommand{\Num}{\refstepcounter{num}%
\textbf{\arabic{section}.\arabic{num}}}

 \newcommand{\Lemma}{\textbf{Lemma~}}
 \newcommand{\Prop}{\textbf{Proposition~}}
 \newcommand{\Cor}{\textbf{Corollary~}}
 \newcommand{\Theorem}{\textbf{Theorem~}}
 \newcommand{\Proof}{\textbf{Proof}}
 \newcommand{\Def}{\textbf{Definition~}}

\date{}
\title{\bf Representations of  solvable  Lie algebras with filtrations}
\author{\bf {A.N.Panov}
\thanks{The paper is supported by  RFBR grants 08-01-00151-a,
09-01-00058-a and by ADTP grant 3341}}
 \maketitle

\section{Introduction. Main definitions}

This paper is devoted to  classification of irreducible
representations of Lie algebras. Complete classification is known
only for Lie algebras of small dimension ~\cite{B1,B2}. This
classification is not  constructible and cannot be generalized to
more  complicated Lie algebras. A more reasonable  approach is to
classify  not irreducible representations but their kernels
(primitive ideals) in the universal enveloping algebra. Started from
 ~\cite{Dix} there were many works in this direction.

An other possibility is to  study some category of
represen\-ta\-tions and irredu\-cible representation. For instance,
the category of represen\-ta\-tions of highest  weight for
semisimple Lie algebras. Irreducible representations in this
category are irreducible factors of Verma modules ~\cite{Dix, BGG}.

In the case of solvable Lie algebras  one well known a family of
irreducible representations  $M_s(f)$, induced from the Vergne
polarization ~\cite{MV}.  In this paper we consider the category of
representations of  Lie algebras with fixed filtrations. We prove
that every irreducible representation in this category have the form
$M_s(f)$  for some  $f\in\gx^*$ (theorem   \ref{IM}).
 That is the mapping  $M_s: f\mapsto M_s(f)$ from  $\gx^*$ to
 $\Irr(\gx,s)$ is surjective. The mapping  $M_s$ is extended to bijection of the factor set
  $\gx^*/\Rx$  to
 $\Irr(\gx,s)$ (theorem \ref{MR} and definition  \ref{DR}).
  Remark that two representations of the form  $M_s(f)$  are equivalent in the category
 of representations of  Lie algebras with filtrations (i.e. $\Rx$-equivalent) if and only if they are
 equivalent in the usual sense (theorem
 \ref{MR}). The classes of $\Rx$-equivalence are described in Theorem  \ref{RR}.

  Using   classes of $\Rx$-equivalence on $\gx^*$ we
   find spectra of the representation induced by irreducible
  representation of subalgebra (theorem \ref{SpInd}),  spectra of restrictions of irreducible representation on subalgebra,
  (theorem \ref{SpRes}) and  spectra of the tensor product of two irreducible repre\-sen\-ta\-tions (theorem \ref{SpT}).

In the last section the connection between the map $M_s:
\gx^*/\Rx\to \Irr(\gx,s)$ and the  Dixmier mapping is investigated.
 The constructed theory is an analog  A.A.Kirillov's orbit
 method (see ~\cite{Dix,K1,K2}) in the category of Lie
algebras with filtration.

 In 1979 this paper was deposed in All-Union Institute for Science and Technical Information (VINITY), see ~\cite{P}.
  Turn to formulation of main definitions of the paper.

1.1. Let  field $K$ be a field of zero characteristic and $\gx$ be a
solvable Lie algebra over $K$. A filtration  $s_\gx$ on $\gx$ is a
chain of ideals

\begin{equation}\label{gs} \gx=\gx_0\supset
\gx_1\supset\ldots\supset\gx_{k-1}\supset\gx_k,
\end{equation}

such that  $\gx_k=\{0\}$ and codimension of $\gx_i$ in $\gx_{i-1}$
as less or equal to 1. A Lie algebra with filtration is a Lie
algebra  $\gx$ with the its fixed filtration  $s_\gx$. To simplify
notations we denote filtration by $s$ and a Lie algebra  $\gx$ with
filtration $s$ by   $(\gx,s)$.

{\bf Remark}. 1) A filtration on $\gx$ always exists, since any
finite dimensional (in particular, adjoint) representation of $\gx$
can be performed in  triangular form   ~\cite[1.3.12]{Dix}.

2) In this definition  we allow that some ideal may  appear several
times. This makes it possible to restrict a  filtration on a
subalgebra (see 1.6) and to  consider a restriction of a
representations of a Lie algebra with filtration  on a subalgebra.

 We say that two filtrations on a Lie algebra are equal, if they
 consists of the same ideals and may differ only by multiplicity of
 their entrance in filtration.

 Let  $(\gx',s')$ and $(\gx'', s'')$  be two Lie algebras with filtrations.
A homomorphism $\phi:(\gx',s')\to(\gx'',s'')$ of Lie algebras with
filtrations is a homomorphism  $\phi:\gx'\to\gx''$ such that
$\phi(\gx'_i)\subset \gx''_i$.

1.2. A chain of ideals $\gx_0\supset \gx_1\supset\ldots\supset\gx_k$
produces the chain of subalgebras  $U(\gx_0)\supset
U(\gx_1)\supset\ldots\supset U(\gx_k)$ of the universal enveloping
algebra  $U(\gx)$. Here we put $U(\gx_k)=K$, since $\gx_k=\{0\}$. We
denote by $U(\gx,s)$ the  universal enveloping algebra  $U(\gx)$
with the fixed filtration.

1.3. Let  $V$  be a linear space over a field $K$ ( infinite
dimensional, in general). A filtration $s_V$ in  $V$ is a chain of
embedded subspaces $$ V=V_0\supset V_1\supset \ldots\supset V_k.
$$
To simplify notations we denote a filtration  by $s$ and a Lie
algebra with filtration $s$ by  $(V,s)$.

{\bf Remark}. Note here that in these case  there is no restrictions
on the codimen\-sions of  $V_i$ in $V_{i-1}$ and on dimensions of
$V_k$.

Assume that for every $i$ we are given a representation  $\rho_i$ of
the Lie algebra $\gx_i$ in linear space  $V_i$, such that the
subspace $V_{i+1}$ is invariant with respect to the restriction of
$\rho_i$ on $\gx_{i+1}$. Then we say that this is a representation
$(\rho, s)$ of the Lie algebra  with filtration $(\gx,s)$ or that
$(V,s)$ is a $U(\gx,s)$-module.

 We say that two representations  $(\rho',s)$ and $(\rho'',s)$ of Lie algebras  $(\gx,s)$ in the spaces
 $(V,s)$ and $(W,s)$ are equivalent, if there exists an isomorphism  $C:V\to W$ such that
  $C(V_i)=W_i$ and $C\rho'(x) =
 \rho'' (x)C$ for every $x\in\gx$.

 1.4. A representation $(\rho,s)$ is  irreducible, if every
$\rho_i$ is an irreducible representation of $\gx_i$ in $V_i$. We
denote the set of equivalence classes of irreducible representations
of $(\gx,s)$  by  $\Irr(\gx,s)$.

1.5. If $(V,s)$ is a linear subspace with filtration and $W$ is a
subspace in  $V$, then the filtration of $V$ induces a filtration of
$W$:
$$W=W_0\supset W_1\supset\ldots\supset W_k,$$
where $W_i= V_i\cap W$.
 The factor space  $V/W$ is also equipped with a natural filtration

$$
(V_0+W)/W \supset (V_1+W)/W \supset\ldots \supset (V_k+W)/W,
$$

that coincides with

$$
V_0/W_0 \supset V_1/W_1 \supset\ldots \supset V_k/W_k.
$$

If  $(V,s)$ is a $U(\gx,s)$-module  and $W$ is a $U(\gx)$-submodule
of  $V$, then we say that  $(W,s)$ is a $U(\gx,s)$-submodule of
$(V,s)$.

A module $T_s$ over $U(\gx,s)$ is  called  a subfactor module of
$(V,s)$, if there exists two  $U(\gx,s)$-submodules  $(W,s)$ and
$(W',s)$ such that  $W\supset W'$ and factor of $(W,s)$ by $(W',s)$
is isomorphic to $T_s$. The set of all irreducible subfactor modules
of  $(V,s)$ is called the spectrum of  $(V,s)$ and is denoted by
$\Spec(V,s)$.

1.6. Let $\hx$ be a Lie subalgebra of  $\gx$. Then filtration of
  $\gx$ produces a filtration
$$\gx_0\cap\hx\supset\gx_1\cap\hx\supset\ldots\supset\gx_k\cap\hx$$
of $\hx$. We shall also denote it by  $s$.

If $(\rho,s)$ is a representation of  $(\gx, s)$, then the
restrictions of   $\rho_i$ on $\hx_i=\gx_i\cap\hx$ in $V_i$ defines
a representation of $(\hx,s)$ in $(V,s)$.

1.7. Let  $\hx$ be a subalgebra of  $\gx$ and  let $(W,s)$ be an
arbitrary $U(\hx,s)$-module. Recall that the induced module
$\ind(W,\gx)$ (see ~\cite[5.1]{Dix}) is defined as linear space
$$U(\gx)\otimes_{U(\hx)}W$$ with the left action of $U(\gx)$.
We consider   a chain of  subspaces in $\ind(W,\gx)$
$$
\ind(W_0,\gx) \supset \ind(W_1,\gx_1) \supset\ldots\supset
\ind(W_k,\gx_k).
$$
Since every  $ \ind(W_i,\gx_i)$  is a $U(\gx_i)$-module,
$\ind(W,\gx,s)$ is a  $U(\gx,s)$-module.

1.8. Let us define the tensor product of representations of Lie
algebras with filtrations. Let  $\gx$ be a Lie algebra with
filtration  $s$ given as   $\gx_0\supset \gx_1\supset\ldots\supset
\gx_k$. We denote the direct sum of two copies of $\gx$ by
$\gx\times\gx$ (instead of the  usual notation $\gx\oplus\gx$). The
Lie algebra $\gx\times\gx$ has a filtration
$$\gx_0\times\gx_0\supset \gx_0\times\gx_1 \supset
\gx_1\times\gx_1\supset\ldots\supset \gx_k\times\gx_k,$$ denoted by
 $s\times s$.

If $V$, $W$ are $U(\gx)$-modules, then we denote by $V\times W$ the
linear space  $V\otimes W$ considered as a
$U(\gx\times\gx)$-module.

 If
$(V,s)$ and $(W,s)$ are  $U(\gx,s)$-modules, then $V\times W$ is a
$U(\gx\times\gx, s\times s)$-module with filtration
$$V_0\times W_0\supset V_0\times W_1\supset V_1\times
W_1\supset \ldots\supset V_k\times W_k.$$  The Lie algebra  $\gx$
embedded into $\gx\times \gx$   diagonally and the filtration
$s\times s$ induces  the filtration $\gx_0\supset \gx_1 \supset
\gx_1 \supset\ldots \supset\gx_k$ on  $\gx$ which coincides with
$s$. The restriction of representation
 $V\times W$ to  $(\gx,s)$ turns  $V\times W$ into
a $U(\gx,s)$-module, which we denote by  $(V,s)\otimes (W,s)$.

1.9. Let $\gx$ be an arbitrary Lie algebra, $\gx^*$  the dual space
of  $\gx$ and $f\in\gx^*$. A subalgebra  $\px$ is called a
polarization of $f$, if $\px$ is a maximal isotropic subspace
 with respect to a skew symmetric form  $f([x,y])$, that is,  $f([\px,\px])=0$ and
$\dim\px=\frac{1}{2}(\dim \gx + \dim\gx^f)$, where $\gx^f$ is the
stabilizer of $f$ in $\gx$. The  restriction of $f$ to its
polarization is a one-dimensional representation (character) of the
polarization.

 Let  $(\gx,s)$ be a Lie algebra with filtration. For any element $f\in\gx^*$   one can  construct
 the polarization,   that  is called the Vergne polarization,
 by the following formula

 $$\pv(f) = \sum_{i=0}^k \gx_i^{f_i},$$
 where $f_i$ is a restriction of $f$ to $\gx^*_i$, and $\gx_i^{f_i}$
 is the stabilizer of $f_i$ in $\gx_i$ (see ~\cite{MV}, \cite[1.12.10]{Dix}).

Denote by  $M(f)$   the module over  $U(\gx)$, induced by the
restriction of  $f$ to $\pv(f)$. The module  $M(f)$ is irreducible
~\cite[6.1.1]{Dix}.

The intersection  $\pv(f)\cap\gx_i$ coincides with the Vergne
 polarization $\pv_i(f_i)$ of $f_i$ in $\gx_i$
~\cite[1.12.10]{Dix}. Therefore

\begin{equation}\label{Ms}
M(f)=M_0(f_0)\supset M_1(f_1)\supset\ldots\supset M_k(f_k),
\end{equation}

where $M_i(f_i)$  is the $U(\gx_i)$-module induced from the
restriction of $f_i$ to $\pv_i(f_i)$. Note that, as $\gx_k=\{0\}$,
then $U(\gx_k)=K$ and $M_k(f_k) = Kl$, where $l=1\otimes 1\in M(f)$.

We denote the module  $M(f)$ with the fixed filtration (\ref{Ms}) by
$M_s(f)$. Note that   $M_s(f)$ is an irreducible $U(\gx,s)$-module.

\section{Irreducible representations of  Lie algebra  $(\gx, s)$}

The main result of this section is Theorem  \ref{IM}. To prove we
need several lemmas.

 Let  $\gx_1$ be an ideal of codimension one in $\gx$. For any   $f_1\in\gx_1^*$ we  introduce notation

$$\gx^{f_1} = \{x\in\gx\mid f_1([x,\gx_1])=0\}.$$

\Lemma\Num\label{Lt}. {\it Let  $f_1\in\gx_1^*$ and $\gx^{f_1} +
\gx_1 = \gx$. Choose $t\in\gx^{f_1}$ and $t\notin\gx_1$.  Let
$\px_1$ be a polarization of  $f_1$ in $\gx_1$ that is invariant
with respect to $\ad_t$. Suppose that  the module  $M_1$, being
equal to $\ind(f\vert_{\px_1},\gx_1)$, is irreducible. Then any
$U(\gx)$-submodule in  $\ind(M_1, \gx)$ has the form $U(\gx)P(t)l$,
where $P(t)$ is some polynomial in  $t$ and  $l = 1\otimes 1$.}

\Proof. There exists the natural filtration of $U(\gx_1)$-submodules
in  the $U(\gx_1)$-module $\ind(M_1, \gx)$:
\begin{equation}\label{Filt}
\{0\} = M_1^{(-\infty)}\subset M_1^{(0)}\subset
M_1^{(1)}\subset\ldots \subset M_1^{(n)}\subset \ldots,
\end{equation}
where $M_1^{(n)} = \oplus_{i\leq n} t^iM_1.$ For any nonzero  $a$ in
$\ind(M,\gx)$, we denote by $\deg(a)$ the smallest  $n$ such that
$a\in M_1^{(n)}$.  If $a=0$, then we put $\deg(a) = -\infty$. We
call  $\deg(a)$ the {\it degree} of   $a$.

{\bf Item 1}. Let  $W$ be a nonzero submodule of  $\ind(M_1, \gx)$
and let $a$ be an element of smallest  degree in  $W$. In this item
we shall prove that $W=U(\gx)a$.

Let  $b\in W$. Using the induction method on $\deg(b)$ we shall show
that $b\in U(\gx)a$.   Case $\deg(b)=-\infty$ (i.e. $b=0$) is
trivial. Suppose that any element in  $W$ of degree smaller than $m$
belongs to  $U(\gx)a$. We shall prove the statement for elements of
degree  $m$. By assumption we have $m\geq n$.

Let  $b=t^mb_0+\ldots+b_m$, $a=t^na_0+\ldots+a_n$, where $a_i,
b_j\in M_1$ and $a_0\ne 0$, ~$b_0\ne 0$. Since $M_1$ is an
irreducible
 $U(\gx_1)$-module, there exists  $u\in U(\gx_1)$, such that $b_0=ua_0$.
 Then the element $c=t^{m-n}ua-b$ has degree less than $m$ and belongs to  $W$.
 By induction assumption,  $c\in U(\gx)a$. Hence
$b\in U(\gx)a$. This proves $W=U(\gx)a$.

{\bf Item 2}. The element $a$ from item  1 can be presented in the
form
$$
a= (t^nv_0+t^{n-1}v_1+\ldots+v_n)l,$$ where $v_i\in U(\gx_1)$. Since
 $M_1$ is an irreducible $U(\gx_1)$-module, there exists  $u_0\in
U(\gx_1)$ such that  $u_0v_0l=l$.  The element $a'$, that is equal
to  $u_0a$, belongs to $W $ and  generates  $W$ as a
$U(\gx_1)$-module (indeed,  any element of degree  $n$ in $W$
generates  $ W$) and have the form
$$ a' = (t^n+t^{n-1}w_1+\ldots+w_n)l,$$
where $w_i\in U(\gx_1)$.

Since the polarization  $\px_1$ is invariant with respect to $\ad_t$
and $f_1([t,\px_1])=0$, we have $$(y-f_1(y))t^nl =
(t-\ad_t)^n(y-f_1(y))l = 0$$ for every $y\in\px_1$. Therefore

$$ (y-f_1(y))a' = (y-f_1(y))(t^n+t^{n-1}w_1+\ldots+w_n)l =
t^{n-1}(y-f_1(y))w_1l + a_1', $$ where $\deg a_1'<n-1$. Since degree
of the element  $(y-f_1(y))a'$ is less than  $n$, we have
$(y-f_1(y))a'=0$ and hence $(y-f_1(y))w_1l=0$ for any $y\in\px_1$.
There exists a unique endomorphism  $\phi$ of the space $M_1$,
commuting with representation of $\gx_1$ in $M_1$, and such that
$\phi(l) =w_1l$. Since  $M_1$ is irreducible, the endomorphism
$\phi$ is scalar  ~\cite[2.6.5]{Dix}. Therefore  $w_1=\al_1\in K$
and $a'=(t^n+t^{n-1}\al_1+t^{n-2}w_2+\ldots+w_n)l$. Arguing as above
one can prove step by step that $w_2,\ldots, w_n$ belong to $K$.
That is $W=U(\gx)P(t)l$ for some polynomial  $P(t)$ over the field
$K$. $\Box$

\Cor\Num\label{CLt}. {\it Let $\gx_1$,~ $M_1$, ~$f_1$ be as in lemma
\ref{Lt}. Any maximal submodule in $\ind(M_1,\gx)$ has the form
$U(\gx)(t-\al)l$ for some $\al\in K$.}

\Lemma\Num\label{Lp}. {\it Let  $\gx$,~$\gx_1$,~$f_1$  be as in the
Lemma
 \ref{Lt}. Denote  by $\pi_1$  projection
$\gx^*\to\gx_1^*$. Let  $\px_1$ be a  polarization of $f_1$ in
$\gx_1$. Suppose that  $\gx^{f_1}+\gx_1=\gx$. Then
\\
1)~ $\gx^f+\gx_1=\gx$ for any $f\in\pi_1^{-1}(f_1)$;\\
2) ~ if $\px$ is a polarization for some  $f\in\pi_1^{-1}(f_1)$ and
$\px\cap\gx_1=\px_1$, then $\px=\gx^{f_1}+\px_1$.}

\Proof.  1) It suffices to prove that if  $x\in\gx^{f_1}$ and
$x\notin \gx_1$, then $x\in \gx^f$. Indeed,
$$f([x,\gx]) =f([x,Kx+\gx_1]) = f([x,\gx_1]) = f_1([x,\gx_1]) = 0.$$
2) Every polarization contains  the stabilizer  $\gx^f$. Hence
$\px=Kx\oplus\px_1$, this proves  2). $\Box$

\Theorem\Num\label{IM}. {\it Every irreducible  $U(\gx,s)$-module
has the form  $M_s(f)$ for some  $f\in\gx^*$.}

\Proof. As we  saw in 1.9, $M_s(f)$ is an irreducible
$U(\gx,s)$-module. We shall prove the theorem using induction on the
 length $k$ of filtration. For  $k=0$ the statement is obvious. To
conclude the proof it suffices to show that the following statement
is true.

{\it Let  $\gx_1$ be the first ideal in filtration (\ref{gs}),
~$f_1\in \gx_1^*$ and  $M_1$ be an irreducible  $U(\gx_1)$-module
induced by the character  $f_1\vert_{\pv_1(f_1)}$ of the Vergne
polarization. Let  $M$  be an irreducible $U(\gx)$-module, that
contains  $M_1$ as an  $U(\gx_1)$-submodule. We require to prove
that  $M=M(f)$ for some $f\in\pi^{-1}(f_1)$.}

The embedding  $M_1$ into $M$ extends to $U(\gx)$-homomorphism
$\Psi: \ind(M_1,\gx)\to M$. Since the module $M$ is irreducible, we
have $\Img(\Psi)=M$.

Consider two cases  a)~ $\gx^{f_1}\subset \gx_1$, ~ b)~
$\gx^{f_1}+ \gx_1 =\gx$.\\
 a)~ $\gx^{f_1}\subset \gx_1$. As  $\gx^f\subset \gx^{f_1}$, we have  $\pv(f)\subset \gx_1$  for
 any $f\in\pi^{-1}(f_1)$.  The Vergne polarization $\pv(f)$ coincides with
 the Vergne polarization of
 $f_1$ in $\gx_1$.  Therefore  $M(f) = \ind(M_1,\gx)$ and
 $\Psi:M(f)\to M$. Both modules  $M(f)$ and  $M$ are irreducible; the Schur Lemma implies
 that  $\Psi$ is  an isomorphism of $M(f)$ onto $M$.\\
  b)~
$\gx^{f_1}+ \gx_1 =\gx$. The polarization  $\px_1=\pv_1(f_1)$
satisfies the conditions of lemma  \ref{Lt} (see ~\cite[1.12.10,
6.1.1]{Dix}). The kernel of homomorphism  $\Psi$ is a maximal
submodule in $\ind(M_1,\gx)$. By Corollary  \ref{CLt}, we obtain
that  the kernel $\Psi$ coincides with $U(\gx)(t-\al)l$ for some
$\al\in K$. Let $f$ be an element of  $\gx^*$ such that  $f(t)=\al$
and $\pi_1(f)=f_1$. By Lemma  \ref{Lp}, $t\in\gx^{f}$  and hence
$\pv(f)= Kt+\pv_1(f_1)$. The module   $M(f)$ coincides with factor
of $\ind(M_1,\gx)$ by $\Ker(\Psi)$ and, therefore, $\Psi$
 isomorphically maps  $M(f)$ onto  $M$. $\Box$

\section{Mapping $M_s$}

In this section we shall answer a question when two modules
$M_s(f')$ and  $M_s(f'')$ are equivalent.

\Lemma\Num\label{Lxf}. {\it There exists a unique nonzero (up to
scalar multiple) element $l$ in $M(f)$ such that
\begin{equation}\label{xf}(x-f(x))l=0.
\end{equation}}
\Proof.  This element exists, since the element   $l=1\otimes 1$
satisfies this property. On the other hand, if $\tilde{l}$ satisfies
(\ref{xf}), then there exists commuting with representation
endomorphism $\phi$ such that  $\phi(l)=\tilde{l}$. As $M(f)$ is an
irreducible  $U(\gx)$-module, the endomorphism  $\phi$ is scalar
~\cite[2.6.5]{Dix}. $\Box$

\Lemma\Num\label{Lyf}. {\it If an element $l$  satisfies (\ref{xf}),
then the equality  $(y-\la)l=0$ implies $y\in\pv(f)$ and $\la=f(y)$,
where $y\in \gx$ и $\la\in K$.}

\Proof. By the Poincar\'e-Birkhoff-Witt theorem (see
~\cite[2.7]{Bur}),  the system
$$\left\{x_1^{k_1}\cdots x_n^{k_n}  ~:~ k_1,\ldots,k_n\in\Zb_+\right\}$$
is a basis of $U(\gx)$ for any basis  $x_1,\ldots, x_n$ of $\gx$.
Let $x_1,\ldots, x_m, x_{m+1},\ldots, x_n$ be a basis of  $\gx$ such
that
 $x_{m+1},\ldots, x_n$ is a basis of  $\pv(f)$. Then $$\left\{x_1^{k_1}\cdots
 x_m^{k_m}l ~:~
k_1,\ldots,k_m\in\Zb_+\right\}$$ is a basis of  $M_s(f)$. $\Box$

\Lemma\Num\label{zf}. {\it Let  $f\in\gx^*$ and  $f_1$ be a
restriction
of  $f$ on the ideal  $\gx_1$ of codimension one from (\ref{gs}). Then \\
a) if  $\gx^f+\gx_1=\gx$, then $M(f)$ is isomorphic to  $M_1(f_1)$
as a
$U(\gx_1)$-module;\\
b) if $\gx^f\subset \gx_1$, then $M(f)$ admits an infinite
filtration, where  each factor  is isomorphic to $M_1(f_1)$.}

\Proof. The statement a) is obvious. Suppose that the assumption of
statement b) holds. Choose $t\in \gx\setminus \gx_1$. As in proof of
Lemma  \ref{Lt},  consider filtration (\ref{Filt}) with $M_1 =
M_1(f_1)$.

Consider the map  $\phi: M_1\to M_1^{(n)}$ such that  $\phi(v) =
t^nv$, ~$v\in M_1$. For any  $u_1\in U(\gx_1)$ we have $u_1\phi(v) =
u_1t^nv = t^nu_1v\bmod M_1^{(n-1)}$. The map  $\phi$ is a
homomorphism of  $U(\gx_1)$-modules. It can de extended  to an
isomorphism  of  $M_1$ to the factor of $M_1^{(n)}$ by
$M_1^{(n-1)}$. $\Box$

\Theorem\Num\label{MR}. {\it The following conditions are equivalent:\\
1)~ $M_s(f')$  and $M_s(f'')$ are equivalent as  $U(\gx)$-modules;\\
2)~ $M_s(f')$  and $M_s(f'')$ are equivalent as  $U(\gx,s)$-modules;\\
3) the Vergne polarizations of $f'$ and $f''$ coincide, and $f'$
equals to
$f''$ under  restriction to a  common Vergne polarization.}\\

\Proof. Implications   $3)\Rightarrow 2) \Rightarrow 1)$ is obvious.
Let  us prove  $1) \Rightarrow 3)$ by induction  on $k$.  For $k=1$
the statement if easy. Suppose that the statement is  proved for
filtrations of length less than $k$. Let us prove for  $k$. Denote
by $s_1$ the restriction of the filtration $s$ to $\gx_1$.

Assume that  $M_s(f')$  and $M_s(f'')$ are equivalent as
$U(\gx)$-modules. Then they are equivalent as  $U(\gx_1)$-modules.
By Lemma \ref{zf},
 $f'$ and $f''$ satisfy simultaneously   either condition  a), or b) of this Lemma,
 and  $U(\gx_1)$-modules $M_{s_1}(f_1')$ and
 $M_{s_1}(f_1'')$  are
 isomorphic, where $f_1'$ and $f_1''$ are restrictions of $f'$ and $f''$ to $\gx_1$.
  Using the the inductive assumption, we obtain:\\
i)~ $\pv_1(f_1') = \pv_1(f_1'')$ (denote it by $\px_1$),\\
ii) ~$f'\vert_{\px_1} = f''\vert_{\px_1}$.

If $f'$ and $f''$ satisfy  condition b) of Lemma \ref{zf}, then
$\px_1$ is their common Vergne polarization; and ii) proves  3).

If  $f'$ and $f''$ satisfy condition  a), then an isomorphism
 $\phi: M_s(f')\to M_s(f'')$ as $U(\gx)$-modules is an isomorphism  $M_{s_1}(f_1')\to M_{s_1}(f_1'')$
 of
$U(\gx_1)$-modules. The conditions  i) and ii) imply that
$M_{s_1}(f_1') = M_{s_1}(f_1'')=\ind(f_1'\vert_{\px_1},\gx_1)$.
Since  $\phi$ is an isomorphism of irreducible module, the operator
$\phi$ is scalar.   Therefore $M_s(f')$ coincides  with $M_s(f'')$
as a $U(\gx)$-module.

Let $l$ be an element of this module, satisfying  (\ref{xf}). By
$M_s(f') = M_s(f'')$, the condition  (\ref{xf}) holds for both cases
$f=f'$,~ $x\in\pv(f')$ and  $f=f''$, $x\in\pv(f'')$. The Lemma
\ref{Lyf} implies that  $\pv(f') = \pv(f'')$ and  $f'$  coincides
with $f''$ under  restriction to a common Vergne polarization.
$\Box$

\Def\Num\label{DR}. Consider the equivalence relation  $\Rx$ on
$\gx^*$ such that  $f'\Rx f''$ if $f'$ and  $f''$ have a common
Vergne polarization and coincide under  restriction to it. We denote
the equivalence class of an arbitrary element $f\in\gx^*$ by
 $\Rx(f)$.

By theorem  \ref{MR},  the correspondence  $M_s:f\mapsto M_s(f)$
extends to bijection $\gx^*/\Rx$ to $\Irr(\gx,s)$. Our next goal is
to describe equivalence classes for  $\Rx$.

\Lemma\Num\label{Lq}. {\it Let  $\gx_1$ be an ideal of codimension
one in
 $\gx$, ~$\px$ (resp. $\px'$) is some polarization  of
$f\in\gx^*$ (resp. $g\in\gx^*$) such that  $\px\cap\gx_1 =
\px'\cap\gx_1=\px_1$, where $\px_1$ is a common polarization of
projections   $f_1$  and  $g_1$ for elements  $f$ and $g$ to
$\gx^*_1$. Let $f\vert_{\px_1} = g\vert_{\px_1}$. Then $\px=\px'$.}

\Proof. The polarizations  $\px$ and $\px'$ can be written in the
form $\px=Kx+\px_1$ and $\px'=Ky+\px_1$ for some  $x,y\in\gx $.
Since $f\vert_{\px_1} = g\vert_{\px_1}$, we have
$$f([y,\px_1]) = g([y,\px_1]) = 0.$$ Therefore, $\px'$ is an isotropic subspace not only for
$g$, but for  $f$ too. Similarly,   $\px$ is an isotropic subspace
for $g$. In particular, the case $\px=\px_1\ne\px'$ as well as
$\px'=\px_1\ne\px$ is not possible.

Hence either $\px=\px'=\px_1$ (that implies the statement of the
lemma), or $\px\ne\px_1$ and $\px'\ne\px_1$. In the second case,
$x,y\notin\gx_1$ and one can choose $x,y$ such that $x-y\in\gx_1$.
Since
$$f_1([x-y,\px_1]) = f([x-y,\px_1]) = f([x,\px_1]) - f([y,\px_1])
=0,$$  the space  $K(x-y)+\px_1$ is isotropic for  $f_1$ in $\gx_1$.
As $\px_1$ is a polarization of $f_1$, we have  $x-y\in\px_1$  and,
finally, $\px=\px'$. $\Box$

\Theorem\Num\label{RR}. {\it The equivalence class  $\Rx(f)$ of an
element  $f\in\gx$ coincides with $\pi^{-1}\pi(f)$, where $\pi$ is
projection  $\gx^*$ to $\pv(f)^*$.}

\Proof. By theorem \ref{MR},  $\Rx(f)\subset \pi^{-1}\pi(f)$. We
have to show that  $\pv(g) = \pv(f)$ for every  $g\in
\pi^{-1}\pi(f)$. Using induction on the length of filtration, we
obtain $\pv_1(f_1) = \pv_1(g_1)$. The  Lemma \ref{Lq} implies
$\pv(f) = \pv(g)$. $\Box$

\section{Spectrum of certain  $(\gx,s)$-modules}

Let  $\hx$ be a subalgebra of Lie algebra  $\gx$. A filtration  $s$
on $\gx$ induces a filtration on   $\hx$ and turns it into  a Lie
algebra $(\hx,s)$ with filtration. Let   $h$ be an element in
$\hx^*$. Denote by $\qv(h)$  the Vergne polarization  of $h$ in
$(\hx,s)$ and by $N_s(h)$  -- the representation
$\ind(h\vert_{\qv(h)},\hx)$. Since  $N_s(h)$ is a representation of
the Lie algebra  $(\hx,s)$, then $\ind(N_s,\gx)$ is a representation
of  $(\gx,s)$ (see 1.7). The goal of this section is to describe
spectra of the induced representation $ \ind(N_s,\gx)$, restriction
$M_s(f)\vert_{(\hx,s)}$ and of representation $M_s(f')\otimes
M_s(f'')$.

\Lemma\Num\label{Zp}.{\it  Let  $Z_s$ be a submodule  of
$\ind(N_s(h),\gx)$;  the factor module by  $Z_s$ is isomorphic to
$M_s(f)$ as a $U(\gx,s)$-module.
Then\\
1) ~$\pv(f)\supset \qv(h)$;\\
 2)~ $f\vert_{\qv(h)} = h\vert_{\qv(h)}$;\\
 3) ~ $Z_s = \sum_{x\in\pv(f)} U(\gx)(x-f(x))l_h$, where $l_h$ is an element
 $1\otimes 1$.}

 \Proof. ~ Use the induction on  the
 filtration length  $k$.  For
 $k=0$  the statement is obvious.
Assume  that the statement is true for length less than  $k$. Let us
 prove it for the length  $k$. Let  $\gx_1$  be an ideal of
  codimension one in filtration
(\ref{gs}). Introduce the notations:  $\hx_1 = \hx\cap\gx_1$, ~$h_1$
(resp. $f_1$) -- projection of $h$ (resp. $f$) to $\hx_1$ (resp.
$\gx_1$), ~$V_s = \ind(N_s(h),\gx)$,~ $s_1$ -- restriction of
filtration $s$ to $\gx_1$, ~ $V_{s_1} = \ind(N_{s_1}(h_1),\gx_1)$,~
$Z_{s_1} = Z_s\cap V_{s_1}$.

By assumption, $V_s/Z_s \cong M_s(f)$. Hence $V_{s_1}/Z_{s_1} \cong
M_{s_1}(f_1)$. By the inductive assumption: \\
1') ~$\pv(f_1)\supset \qv(h_1)$;\\
2')~ $f_1\vert_{\qv(h_1)} = h_1\vert_{\qv(h_1)}$;\\
3') ~ $Z_{s_1} = \sum_{x\in\pv(f_1)} U(\gx)(x-f_1(x))l_{h_1}$, where
$l_{h_1} =l_h=1\otimes 1$. Consider two cases: a)
~$\qv(h)+\gx_1=\gx$ и b) ~$\qv(h)\subset
\gx_1$. \\
a) ~$\qv(h)+\gx_1=\gx$. In this case $V_s=V_{s_1}$ and, therefore,
$Z_s=Z_{s_1}$, that is  $Z_{s_1}$ is a $U(\gx)$-submodule.

Let $x\in\pv_1(f_1)$. It follows from  3') that $(x-f_1(x))l_h\in
Z_{s_1}$. Since  $Z_{s_1}$ is a  $U(\gx)$-module, $y(x-f_1(x))l_h\in
Z_{s_1}$ for every $y\in\gx$. Choose $y\in\qv(h)$ and $y\notin
\gx_1$. The following equality
$$ y(x-f_1(x))l_h = (x-f_1(x))yl_h + [y,x]l_h = h(y)(x-f_1(x))l_h +
[y,x]l_h$$ implies  $[y,x]l_h\in Z_{s_1}$.   After factorization by
$Z_{s_1}$  we obtain that $[y,x]l =0$ holds in $M_{s_1}(f_1)$,where
$l=1\otimes 1$.
        By lemma  \ref{Lyf}, we have got $[y,\pv_1(f_1)]\subset \pv_1(f_1)$
        and
$f_1([y,\pv_1(f_1)]) = 0$.

So, the subalgebra  $\px$, defined as $Ky+\pv(f_1)$, is a
polarization for every  $g\in\pi_1^{-1}(f_1)$, where $\pi_1$ is
projection of $\gx^*$ to $\gx_1^*$. Choose $g\in\pi_1^{-1}(f_1)$
such that $g(y) = h(y)$.
Using inductive assumption, we obtain \\
1'') ~$\px\supset \qv(h)$;\\
 2'')~ $g\vert_{\qv(h)} = h\vert_{\qv(h)}$.\\
Since  $y\in\qv(h)$, we have $(y-g(y))l_h = (y-h(y))l_h = 0$. Recall
that in our case  $Z_s = Z_{s_1}$. Applying 3'), we obtain
 $$ \mbox{3'')} ~~ Z_{s_1} = \sum_{x\in\pv_1(f_1)} U(\gx)(x-f_1(x))l_{h} = \sum_{x\in\px} U(\gx)(x-g(x))l_{h}.$$

Since $\px\cap\gx_1=\pv_1(f_1)$ and  $\pv(g)\cap\gx_1=\pv_1(f_1)$,
we have $\px=\pv(g)$ (see Lemma  \ref{Lq}). As $V_s =
\ind(N_s(h),\gx) = \ind(h\vert_{\qv(h)},\gx)$ and  $Z_s$ has the
form  3''), we conclude $V_s/Z_s\cong M_s(g)$. By assumption, we
have $V_s/Z_s\cong M_s(f)$. Therefore, $M_s(f)\cong M_s(g)$. Using
theorem  \ref{MR}, we obtain $\pv(g)=\pv(f)=\px$ and $f\vert_\px =
g\vert_\px$. Substituting this equalities in 1''), 2''), 3''), we
get 1), 2), 3). This proves lemma
in the case  a).\\
b) ~$\qv(h)\subset \gx_1$.  In this case, $V_s = \ind(V_{s_1},\gx)$,
~~ $Z_s\supset\ind(Z_{s_1},\gx)$ and
$$V_s/\ind(Z_{s_1},\gx) = \ind(V_{s_1},\gx)/\ind(Z_{s_1},\gx) = \ind(M_{s_1}(f_1),\gx).$$

Denote by  $\overline{Z}_s$ the image of $Z_s$ in factor of  $V_s$
with respect to
$\ind(Z_{s_1},\gx)$. Consider two cases.\\
b1)~ $ \gx^{f_1}\subset\gx_1$. In this case, $\pv(f)=\pv_1(f_1)$
(this proves  1) and 2)) and  $U(\gx)$-module
$\ind(M_{s_1}(f_1),\gx)$ is irreducible. Hence,  $\overline{Z}_s=0$
and, therefore,
$$ Z_s = \ind(Z_{s_1},\gx) = \sum_{x\in\pv(f)=\pv_1(f_1)}
U(\gx)(x-f_1(x))l_{h};$$ this proves 3).\\
b2)~ $\gx^{f_1}+\gx_1 = \gx$. In this case, the Vergne polarization
$\pv(f)$ has the form  $Ky\oplus\pv_1(f_1)$, where
$y\in\pv(f)\setminus \gx_1$ (see  lemma \ref{Lq}). By corollary
\ref{CLt}, we get  $\overline{Z}_s = U(\gx)(y-f(y))l_h$, that
implies  the statements  1), 2) and 3). $\Box$

\Cor\Num\label{CZp}. {\it The $U(\gx,s)$-module  $M_s(f)$ is
isomorphic
to factor of $\ind(N_s(h),\gx)$ if and only if  \\
1)~ $\pv(f)\subset \qv(h)$ and 2) ~$f\vert_{\qv(h)} =
h\vert_{\qv(h)}$.}

\Proof. It is immediate from Lemma \ref{Zp}. $\Box$

\Theorem\Num\label{SpInd}. {\it Let $\rx(h)$ be the equivalence
class of
 $h\in\hx^*$,~
$\pi$ be the projection  $\gx^*$ на $\hx^*$.  Then
$$\Spec(\ind(N_s(h),\gx)) = \left\{M_s(f)~\vert~ \Rx(f)\subset
\pi^{-1}(\rx(h))\right\}.$$ Here, as above, $\Rx(f)$  is the
equivalence class of $f\in\gx^*$.}

 \Proof. The $U(\gx,s)$-module
$V_s(h)$ equals to   $\ind(N_s(h),\gx)$ and has the filtration
$$V_0(f_0)\supset V_1(f_1)\supset \ldots\supset V_k(f_k),$$
here  $V_k(h_k)=kl$, ~ $l=1\otimes 1$.

Suppose that $M_s(f)\in\Spec(V_s(h))$. This means that there exist
$U(\gx,s)$-submodules $W_s$ and $W_s'$ such that  $V_s\supset
W_s\supset W'_s$ and $W_s/W'_s$ is isomorphic to  $M_s(f)$. This is
equivalent to the existence of chains $W_0\supset W_1\supset
\ldots\supset W_k$ and $W'_0\supset W'_1\supset \ldots\supset W'_k$
such that $W_0$,~ $W'_0$  are $U(\gx)$-modules, ~ $W_i= W_0\cap
V_i(h_i)$,~ $W'_i= W'_0\cap V_i(h_i)$,~ $V_i(h_i)\supset W_i\supset
W'_i$ and $W_i/W'_i$ is isomorphic to $M_i(f_i)$ as
$U(\gx_i)$-module. For $i=k$ we have $V_k(h_k)\supset W_k\supset
W'_k$ and $W_k/W'_k$ is isomorphic to $M_i(f_i)$ as a linear space
over  $K$. Since $\dim V_k(h_k) = \dim M_i(f_i) = 1$, we conclude
that  $W'_k=\{0\}$ and $W_k=V_k(h_k)=Kl$. Hence $l\in W_s$. By
$V_s(h)=U(\gx)l$,  $W_s= V_s$ and any  $U(\gx,s)$-subfactor module
is a factor module.

This shows that     $ M_s(f)\in \Spec\left(V_s(h)\right)$ if and
only if  $f$ and  $h$ satisfy the conditions 1) and 2) of Corollary
 \ref{CZp}.

We have
$$
\Rx(f) = \left\{ \xi\in\gx^* \mid \xi\vert_{\pv(f)} =
f\vert_{\pv(f)}\right\},$$
$$
\pi^{-1}(\rx(h)) = \left\{ \eta\in\gx^*\mid \eta\vert_{\qv(h)} =
h\vert_{\qv(h)}\right\}$$ and conditions  1) and  2) of Corollary
\ref{CZp} are equivalent to  $\Rx(f)\subset \pi^{-1}(\rx(h))$.
$\Box$

\Theorem\Num\label{SpRes}.{\it Let $\rx(h)$, ~ $\Rx(f)$,~ $\pi$ be
as in previous theorem. Then
$$ \Spec\left( \left. M_s(f)\right\vert_{(\hx,s)}\right) = \left\{ N_s(h)~ \vert~ \rx(h)
\subset \pi(\Rx(f))\right\}.$$}

\Proof. Suppose that  $N_s(h)\in\Spec\left.
M_s(f)\right\vert_{(\hx,s)}$. This means that there exist
$U(\hx,s)$-submodules  $C_s$ and $C_s'$ such that  $C_s'\supset C_s$
and $C_s/C_s'$ is  isomorphic to  $N_s(h)$. Arguing as in the proof
of previous theorem, we obtain  $C'_k\supset C_k\supset M_k(f_k)$
and $C_k/C_k'\cong
 N _k(h_k)$. Since $\dim M_k(f_k) =\dim N_k(h_k)=1$, we have
 $C_k'=\{0\}$ and $C_k = M_k(f_k)=Kl$, where $l=1\otimes 1$. Hence
 $l\in C_k$ and  $U(\hx)l\subset C_s$. This implies that
 $N _s(h)$ is a factor of  $U(\hx,s)$-module $U(\hx)l$, that is
 isomorphic to  $\ind(f\vert_{\pv(f)\cap \hx}, \hx)$. Applying
 Corollary
\ref{CZp}, we obtain that $N_s(h)$ is contained in
$\Spec\left(\left.M_s(f)\right\vert_{(\hx,s)}\right)$ if and only if

\begin{equation}\label{Rab}
\left\{\begin{array}{l}
 \qv(h)\supset \pv(f)\cap\hx,\\  h\vert_{\pv(f)\cap\hx} =
f\vert_{\pv(f)\cap\hx}.\end{array}\right.
\end{equation}
Since
$$
\rx(h) = \left\{ \la\in\hx^* \mid \la\vert_{\qv(h)} =
h\vert_{\qv(h)}\right\},$$
$$
\pi(\Rx(f)) = \left\{ \mu\in\hx^*\mid \mu\vert_{\pv(f)\cap \hx} =
f\vert_{\pv(f)\cap\hx}\right\},$$ tha conditions  (\ref{Rab}) are
equivalent to  $\rx(h)\subset \pi(\Rx(f))$.

 $\Box$

Let $f', f''\in\gx^*$. The Vergne polarization of $f'\times f''\in
(\gx\times\gx)^*$  with respect to filtration $s\times s $ (see 1.8)
 coincides with   $\pv(f')\times\pv(f'')$ and $$
M_s(f')\times M_s(f'') = M_{s\times s}(f'\times f'').$$

\Theorem\Num\label{SpT}. {\it Let $f',f''\in\gx^*$. we claim that
$$\Spec~ M_s(f')\otimes M_s(f'') = \left\{M_s(g) \mid
\Rx(g)\subset \Rx(f')+\Rx(f'')\right\}.$$}

\Proof. By definition, $M_s(f')\otimes M_s(f'')$ is a restriction of
 $ M_s(f')\times M_s(f'')$ to $(\gx,s)$. By Theorem
\ref{SpRes},
$$ \Spec~ M_s(f')\otimes M_s(f'') = \left\{M_s(g) \mid
\Rx(g)\subset \pi\left(\Rx(f'\times f'')\right) \right\},$$ where
 $\pi$ is  projection  $(\gx\times\gx)^*$ to $\gx^*$, ~ $\Rx(f'\times
 f'')$ is the equivalence class  of $f'\times f''$ in $(\gx\times\gx)^*$.
 Easy to see that  $$\pi\left(\Rx(f'\times f'')\right) = \pi
 \left(\Rx(f')\times \Rx(f'')\right) = \Rx(f') + \Rx(f''). \Box $$

 \section{Connection with Dixmier mapping. Examples}

The orbit method of A.A.Kirillov  describes irreducible
representations of solvable Lie groups in Hilbert spaces. It has an
algebraic analogue. Let  $\gx$ be  a solvable Lie algebra over
algebraically closed field $K$ of zero characteristic and  let $\Ac$
be an adjoint algebraic group for $\gx$. Let $f$ be an  element of
$\gx^*$  and $\px$ be a polarization of  $f$. We denote by
$\theta_\px$ the character of Lie subalgebra  $\px$ given as
$\frac{1}{2}\tr\ad_{\gx/\px}$.
 The ideal  $I(f)$ defined as a  kernel of the twisted induced representation
 $\ind\,\tilde{}\,(f\vert_\px,\gx) =\ind(f - \theta_\px\vert_\px,\gx)$
 in the universal enveloping algebra
 is primitive  and does not depend on choice of polarization  $\px$. The
 correspondence
  $f\mapsto I(f)$   can be extended to  bijection
   of the orbit space  $\gx^*/\Ac$ to that set  $\Prim\,U(\gx)$ of primitive ideals
   of
   $U(\gx)$. This bijection if denoted by   $\Dix_\gx$
   and is called the Dixmier map; see  ~\cite[Глава 6]{Dix}.

\Theorem\Num. {\it 1) $\Rx(f)\subset\Ac(f)$ for  every  $f\in\gx^*$.
Denote by $i: \gx^*/\Rx\to \gx^*/\Ac$ the embedding of
 $\Rx(f)$ into $\Ac(f)$.\\
2)  Denote by  $\KKer: \Irr(\gx,s)\to\Prim ~U(\gx)$ the map sending
the irreducible representation  $\ind (f\vert_\px,\gx)$ to the
kernel of twisted  representation
$\ind\,\tilde{}\,(f\vert_\px,\gx)$. \\ Then  the following diagram
is commutative:
   $$
     \begin{CD}
             \gx^*/\Rx @>i>> \gx^*/\Ac\\
             @VV M_sV @VV \Dix_\gx V  \\
             \Irr(\gx,s) @> \KKer >> \Prim\,U(\gx)
     \end{CD}
   $$
}

   \Proof. 1)  For any  $g\in\Rx(f)$  we have   $\pv(g)=\pv(f)=\px$
   and $g\vert_\px = f\vert_\px$ (theorem \ref{RR}). Hence the
   map   $g\mapsto
   I(g) = \Ker\left(\ind\,\tilde{}\,(f\vert_\px,\gx)\right)$
   is constant on $\Rx(f)$. Therefore  $\Rx(f)\subset \Ac(f)$.

   The statement of item  2) is a corollary of definition of the Dixmier mapping.  $\Box$

{\bf Example 1}.  Let  $\gx=  \spn <x,y\mid [x,y]=y>$. Fix the
filtration  $s$ of the form  $\gx\supset Ky\supset\{0\}$.

For $f\in\gx^*$  denote  $f(y)=y^0$. If $y^0\ne 0$, then $\Rx(f) =
\{g\in\gx^* \mid g(y)=y^0\}$; and if  $y^0 = 0$, then $\Rx(f) =
\{f\}$.

If  $y^0\ne 0$, then the corresponding irreducible representation of
 $(\gx,s)$ is induced  by the character  $y\mapsto y^0$
of  subalgebra  $\px=Ky$. In the case  $y_0=0$ linear form  $f$ is a
character and defines one-dimensional  representation of $(\gx,s)$.
Note that a classification of all irreducible representations of Lie
algebra $\gx$ are given in  ~\cite{B2}.
  \\
  {\bf Example 2}.  Let  $\gx$ be the Heisenberg algebra  $ \spn <x,y, z\mid [x,y]=z>$.
  Fix the filtration  $s$ of the form  $\gx\supset Ky\oplus Kz\supset Kz\supset\{0\}$.

Let $f\in\gx^*$,  denote   $f(y)=y^0$ and $f(z)=z^0$. If  $z^0\ne
0$, then $\Rx(f) = \{g\in\gx^* \mid g(z) = y^0,~~ g(z)=z^0\}$; and
if $z^0 = 0$, then $\Rx(f) = \{f\}$.

If $z^0\ne 0$, then the corresponding irreducible representation of
$(\gx,s)$ is induced by the  character  $y, z\mapsto y^0, z^0$ of
the subalgebra  $\px=Ky\oplus Kz$. In the case  $z_0=0$ the linear
form $f$ is a character. It  defines one-dimensional representation
of $(\gx,s)$.

Concerning a  classification of all irreducible representations of
the Heisenberg algebra, it is known that every irreducible is either
one-dimensional (i.e. coincides with  $f$ with $z^0=0$), or is an
irreducible representation  of the Weyl algebra  $A_1 = <p,q\mid
[p,q]=1 >$ (see ~\cite{B1}).

Note that in each of these  examples irreducible representations of
$(\gx,s)$ cover a small part of all irreducible representations of
$\gx$ (see \cite[Prop. 6.1]{B2}).

A.N.Panov)\\
Samara State University,
\\
{\it E-mail}: apanov@list.ru


\begin{thebibliography}{100}



\bibitem{B1}
 R.E. Block,
 "The irreducible representations of the Weyl algebra $A_1$"{}, {\it Lecture
Notes in Math.}, {\bf 740}(1979)  Springer-Verlag, Berlin/New York,  69-79.

\bibitem{B2}
 R.E. Block,   "The irredicible representations of  the
 Lie algebra  $\mathrm{sl}(2)$ and of the Weyl algebra"{},
{\it Advances in Math.}, {\bf 39}(1981), 69-110.

 \bibitem{Dix}
 J.Dixmier, Alg\`ebres enveloppantes, Gauthier-Villars, Paris, 1974.

\bibitem{BGG}
I.N.Berbstein, I.M.Gelfand, S.I.Gelfand  "Structure of
representations, generated by highest weight vectors"{}, {\it
Funct.analys i prilozh.}, {\bf 5}(1971), 1-9 [in russian].




\bibitem{MV}
M. Vergne, "Construction de sous-alg\`ebres subordonn\'ees a un
\'el\'ement du dual d'une alg\`ebre de Lie r\'esoluble"{}, {\it C.R.
Acad. Sci. Paris(A)}, {\bf 270}(1970), 173-175, 704-707.

\bibitem{K1}
A.A.Kirillov, " Unitary representations of nilpotent Lie groups"{},
{\it Uspehi Math. Sci.}, {\bf 17}(1962), 57-110 [in russian].

\bibitem{K2}
A.A.Kirillov, Lectures on the orbit method, Graduate Studies in
Math.,
 64(2004), Providence, RI: AMS.

 \bibitem{P}
 A.N.Panov, "Representations of Lie algebras with filtration",
 VINITI, no. 3012-79, 1979 [in russian].

\bibitem{Bur}
N.Bourbaki, Groupes et algebras de Lie (chapitre I-III), Hermann,
Paris, 1972.


\end{thebibliography}
\end{document}